\documentclass[12pt,a4paper]{article}
\usepackage{amssymb,amsmath,amsthm}

\def\ex{{\rm ex}}
\def\eps{{\varepsilon}}

\def\BB{{\mathcal B}}
\def\CC{{\mathcal C}}

\def\FF{{\mathcal F}}
\def\GG{{\mathcal G}}

\def\HH{{\mathcal H}}
\def\KK{{\mathcal K}}

\newcommand\cF{{\mathcal F}}
\newcommand\cG{{\mathcal G}}
\newcommand\cH{{\mathcal H}}

\theoremstyle{plain}
\newtheorem{theorem}{Theorem}[section]
\newtheorem{lemma}[theorem]{Lemma}
\newtheorem{corollary}[theorem]{Corollary}
\newtheorem{conjecture}[theorem]{Conjecture}

\theoremstyle{definition}

\newtheorem{prob}[theorem]{Problem}

\newcommand\lref[1]{Lemma~\ref{lem:#1}}
\newcommand\tref[1]{Theorem~\ref{thm:#1}}
\newcommand\cref[1]{Corollary~\ref{cor:#1}}

\newcommand\sref[1]{Section~\ref{sec:#1}}

\begin{document}

\title{Large $B_d$-free and union-free subfamilies}

\author{J\'anos Bar\'at\thanks{Department of Computer Science and Systems Technology,
 University of Pannonia, Egyetem u. 10, 8200 Veszpr\'em, Hungary.  Research is supported by OTKA Grant PD-75837 and
 the J\'anos Bolyai Research Scholarship of the Hungarian Academy of Sciences.} \and
Zolt\'an F\"uredi\thanks{Alfr\'ed R\'enyi Institute of Mathematics, P.O.B. 127, Budapest H-1364, Hungary and
 Department of Mathematics, University of Illinois at Urbana-Champaign. Email: z-furedi@illinois.edu.
 Research supported in part by the Hungarian National Science Foundation
 OTKA, and by the National Science Foundation under grant NFS DMS 09-01276 ARRA.} \and
Ida Kantor\thanks{Department of Mathematics, University of Illinois at Urbana-Champaign, and
 Department of Applied Mathematics, Charles University, Malostransk\'{e} n\'{a}m. 25, 11800 Praha 1, Czech Republic.
 E-mail: ida@kam.mff.cuni.cz.} \and
Younjin Kim\thanks{Department of Mathematics,
 University of Illinois at Urbana-Champaign, 1409 W. Green Street, Urbana, IL 61801. Email: ykim36@illinois.edu.} \and
Bal\'azs Patk\'os\thanks{Alfr\'ed R\'enyi Institute of Mathematics, P.O.B. 127, Budapest H-1364, Hungary.
 Email: patkos@renyi.hu. Research
supported by Hungarian National Scientific Fund, grant number: OTKA K-69062 and PD-83586. }
}
\maketitle

\begin{abstract}
For a property $\Gamma$ and a family of sets  $\cF$ , let $f(\cF,\Gamma)$ be the size
of the largest subfamily of $\cF$ having property $\Gamma$.
For a positive integer $m$, let $f(m,\Gamma)$ be the minimum of $f(\cF,\Gamma)$ over all
families of size $m$. A family $\cF$ is said to be $B_d$-free if it has no subfamily $\cF'=\{F_I: I \subseteq [d]\}$ of $2^d$ distinct
sets such that for every $I,J \subseteq [d]$, both $F_I \cup F_J=F_{I \cup J}$ and $F_I \cap F_J = F_{I \cap J}$ hold.
A family $\cF$ is $a$-union free if $F_1\cup \dots F_a \neq F_{a+1}$ whenever $F_1,\dots,F_{a+1}$ are distinct sets in $\FF$.
We verify a conjecture of Erd\H os and Shelah that $f(m, B_2\text{\rm -free})=\Theta(m^{2/3})$.
We also obtain lower and upper bounds for $f(m, B_d\text{\rm -free})$ and $f(m,a\text{\rm -union free})$.
\end{abstract}

\section{Introduction, results}

Moser proposed the following problem: Let $A_1,A_2\dots,A_m$ be a collection of $m$ sets.
A subfamily $A_{i_1},A_{i_2}\dots,A_{i_r}$ is \textit{union-free} if
$A_{i_{j_1}}\cup A_{i_{j_2}}\neq A_{i_{j_3}}$ for every triple of distinct sets
$A_{j_1},A_{j_2},A_{j_3}$ with $1\leq j_1\leq r$,
 $1\leq j_2\leq r$, and $1\leq j_3\leq r$.
Erd\H os and Koml\'os ~\cite{EK1969} considered the following problem of Moser:
what is the size of the largest union-free subfamily $A_{i_1},\dots,A_{i_r}$?

Put $f(m)=\min r$, where the minimum is taken over all families of $m$
distinct sets.
As mentioned in \cite{EK1969}, Riddel pointed out that $f(m)>c\sqrt m$.
Erd\H os and Koml\'os~\cite{EK1969}
showed $\sqrt{m}\leq f(m)\leq 2\sqrt 2\sqrt{m}$.
Kleitman proved $\sqrt{2m}-1<f(m)$; Erd\H os and Shelah~\cite{ES1972} obtained
\begin{equation}\label{ES_upper}
  f(m) <2\sqrt{m}+1.
  \end{equation}
The latter two conjectured $f(m)=(2+o(1))\sqrt{m}$.

We define $f(\cF, \Gamma)$  as the size of the largest subfamily of $\cF$
having property $\Gamma$,
$$
  f(\cF, \Gamma):=\max \{  |\cF'|: \cF'\subseteq \cF, \enskip  \cF'{\rm
\enskip  has \enskip property \enskip}\Gamma\}.
  $$
In this context, $f(E(K_r^n), \cH\text{\rm -free})$ is the Tur\'an number $\textrm{ex}_r(n, \cH)$.
Let $f(m, \Gamma)= \min \{ f(\cF, \Gamma): |\cF|=m \}$.
Generalizing the union-free property,
a family $\cF$ is  \textit{$a$-union free} if there are no distinct sets $F_1, F_2\dots, F_{a+1}$
satisfying $F_1\cup F_2 \cup \dots \cup F_a=F_{a+1}$.

Erd\H os and Shelah~\cite{ES1972} also considered $\Gamma$ to be the property that no four
 distinct sets satisfy $F_1\cup F_2=F_3$ and $F_1\cap F_2=F_4$.
Such families are called $B_2$-{\em free}.
Erd\H os and Shelah~\cite{ES1972} gave an example showing
 $f(m, B_2\text{\rm -free})\leq(3/2)m^{2/3}$ and
 they also conjectured  $f(m, B_2\text{\rm -free})>c_2 m^{2/3}$.

A family $\BB$ of $2^d$ distinct sets is forming a Boolean algebra of dimension $d$
  if the sets can be indexed with the subsets of $[d]=\{1,2,\dots,d\}$ so that
$B_I \cap B_J=B_{I \cap J}$ and $B_I \cup B_J=B_{I \cup J}$ hold for any $I,J \subseteq [d]$.
If $\FF$ does not contain any subfamily forming a Boolean algebra of dimension $d$, then it is called
 $B_d$-{\em free}, or we say that $\FF$ {\it avoids} any Boolean algebra of dimension $d$. A result by Gunderson, R\"odl, and Sidorenko \cite{GRS} states that
$f(2^{[n]}, B_d\text{\rm -free})=\Theta(2^n/n^{2^{-d}})$.
In \sref{Bd}, we prove the aforementioned conjecture by Erd\H os and Shelah in the following
 more general form.

\begin{theorem}
\label{thm:bdf}
For any integer $d$, $d \ge 2$, there exist constants $c_d, c_d'>0$,  and exponents
$$
     e_d:=\frac{2^d-\lceil \log_2 (d+2) \rceil}{2^d-1}, \quad e_d':=\frac{2^d-2}{2^d-1}
  $$
 such that
$$c_dm^{e_d}  \le f(m, B_d\text{\rm -free}) \le c_d'm^{e_d'}.$$
In particular,  
\begin{equation}\label{eq:2}
 (3\cdot 2^{-7/3}+o(1))m^{2/3} \le f(m,B_2\text{\rm -free}) \le \frac{3}{2}m^{2/3}.
  \end{equation}
\end{theorem}

In \sref{Uf} we consider $a$-union free families.
We generalize the construction giving (\ref{ES_upper})  
  and prove the following

\begin{theorem}
\label{thm:afree}
For any integer $a$, $a \ge 2$,
\begin{equation}\label{eq:3}
\sqrt{2m}-\frac{1}{2} \le f(m,a\text{\rm -union~free}) \le 4a+ 4a^{1/4}\sqrt{m}.
  \end{equation}
\end{theorem}

Since the first version of this manuscript, Fox,  Lee, and Sudakov~\cite{FLS} verified the present authors'
 conjecture (see later as Problem~\ref{prob_limit}) and proved a matching lower bound showing that
 $f(m,a\text{\rm -union~free}) \ge \max\{ a, \frac{1}{3}\sqrt[4]{a}\sqrt{m})\}$.

\section{Subfamilies avoiding Boolean algebras of dimension $d$}
\label{sec:Bd}

In this section we prove the lower bounds in \tref{bdf}
 by a probabilistic argument applying the first moment method.

Suppose that $\BB=\{B_I: I\subseteq [d]\}$ is forming a Boolean algebra  of dimension $d$.
Thus we have nonempty, pairwise disjoint sets, $A_0, A_1, \dots, A_d$, called {\it atoms},  such that
 $B_I=A_0\cup \{ A_i: i \in I\} $.
 A subfamily $\CC\subseteq \BB$  {\it determines} the Boolean algebra $\BB$ if
  every member of $\BB$ can be obtained as a Boolean expression
 (using unions, intersections, differences, but not complements)  of some sets of $\CC$.
Obviously, the $d$ sets of the form  $\{ A_0\cup A_i: i\in [d]\}$ determine $\BB$.
Much more is true.

\begin{lemma}
\label{lem:numbd}
Suppose that the sets of $\BB$ are forming a Boolean algebra  of dimension $d$.
Then there exists a subfamily $\CC\subseteq \BB$ of size  $\lceil \log_2(d+2) \rceil$  and determining  $\BB$.
Moreover, there is no subfamily of smaller size with the same property.
\end{lemma}

\begin{proof}
Let $k:=\lceil \log_2(d+2) \rceil$.
We define an appropriate $\CC$ of size $k$ by considering a  standard construction used for
 non-adaptive binary search.
Namely,  write each integer $i\in [d]$ in base 2, $i=\sum_{1\leq j \leq k} \eps_{i,j}2^{j-1}$ and define
  $C_j=A_0\cup \{ A_i: \eps_{i,j}=1\}$, $j=1,2, \dots, k$.

On the other hand, any $\CC$ determines at most $2^{|\CC|}-1$ nonempty atoms,  we obtain
  $2^{|\CC|}-1\geq d+1$.
\end{proof}

 \begin{corollary}
	\label{cor:sizebd}
Given any family $\cF=\{F_1,F_2,\dots, F_m\}$  of $m$ sets, $\cF$ contains at most $\binom{m}{\lceil\log(d+2)\rceil}$ subfamilies
 forming a Boolean algebra of dimension $d$.
	\end{corollary}

\lref{numbd} gives the right order of magnitude on the number of possible subfamilies
 forming a Boolean algebra of dimension $d$  contained in a family of $m$ sets, as shown by the family
 $\cF=2^{[n]}$, where $m=2^n$.

\begin{proof}[Proof of the lower bound in {\rm \bf \tref{bdf}}.]
Let $\cF=\{F_1,F_2,\dots, F_m\}$ be any family of $m$ sets.
Let us consider a random subfamily $\cF'$, that is, we select every set in $\cF$ independently with probability $p$.
Let $X$ be the random variable denoting the number of sets in $\cF'$, and let $Y$ be the random variable denoting the
number of subfamilies in $\cF'$ forming a Boolean algebra of dimension $d$.
By \cref{sizebd},
$$\mathbb{E}(X-Y) \ge mp-p^{2^d}\binom{m}{\lceil\log_2(d+2)\rceil}.$$
If we remove a set from each subfamily in $\cF'$ forming a Boolean algebra of dimension $d$, then
 we obtain a $B_d$-free subfamily $\cF''$ of size at least $X-Y$.
Substituting $p=m^{e_d}$ where $e_d=\frac{\lceil \log (d+2) \rceil -1}{2^d-1}$ yields the lower bound.
To get a better constant in the case $d=2$, put $p=2^{-1/3}m^{-1/3}$.
\end{proof}

One might try to improve the constant of the lower bound by improving \lref{numbd} for
  families without large chains and antichains.
However,  the construction of Erd\H os and Shelah shows, one cannot hope for anything better
  than $(\frac{1}{2}+o(1))\binom{m}{2}$, which would improve the constant of the lower bound in (\ref{eq:2})
 only to $3/4$.

\section{Upper bound using Tur\'an theory}\label{s:turan}

In this section we prove the upper bounds in \tref{bdf}
 by generalizing the ideas of Erd\H os and Shelah \cite{ES1972}.

Let $\KK(a_1, \dots, a_d)$ denote the
 complete, $d$-partite hypergraph with parts of sizes $a_1, \dots, a_d$,
 i.e., $V(\KK):=X_1\cup \dots\cup X_d$ where $X_1, \dots, X_d$ are pairwise disjoint sets
 with $|X_i|=a_i$, and  $E(\KK):=\{ E: |E|=d, \, |X_i\cap E|=1$ for all $i\in [d]\}$.
For short we use  $\KK_d^{(k)}$
 for $\KK(k,k^2,\dots,k^{2^{d-1}})$ and
 $K_{d\ast 2}$ for $\KK(2,\dots, 2)$.
The (generalized) {\it Tur\'an number} of the $d$-uniform hypergraph $\HH$ with respect
 to the other hypergraph $\GG$, denoted by $\ex(\GG,\HH)$,  is the size of the largest
 $\HH$-free subhypergraph of $\GG$.

\begin{theorem}\label{K22free}
For $k, d\geq 2$, $\ex(\KK_d^{(k)}, K_{d\ast 2})< \left( 2-\frac{1}{2^{d-1}} \right) k^{2^d -2}$.
\end{theorem}

\begin{proof}
 We proceed by induction on $d$.
Let $d=2$, and let $H$ be a $K_{2,2}$-free subgraph of $K_{k,k^2}$.
Let $v_1,v_2,\dots,v_{k^2}$ be the vertices of the larger part of $K_{k,k^2}$, and $d_i:=\deg_H(v_i)$.
Each pair of vertices in the smaller part of $K_{k,k^2}$ has at most one common neighbor in $H$.
Therefore,
$\sum \binom{d_i}{2} \leq \binom{k}{2}$. It yields
\[
 |E(H)|= \sum_{i=1}^{k^2} {d_i} \leq
\sum_{i=1}^{k^2} \left( \binom{d_i}{2}+1\right)
\leq \binom{k}{2} +k^2.
\]

Fix $d$, $d>2$, and a $K_{d\ast 2}$-free subhypergraph $\HH$ of $\KK_d^{(k)}$.
Let $v_i$ $1\leq i\leq k^{2^{d-1}}$ be the vertices of the largest part of $\KK_d^{(k)}$, and $d_i:=\deg_{\HH}(v_i)$.
Let $\HH_i$ be the $(d-1)$-uniform $(d-1)$-partite hypergraph, which we get by taking the set of edges of $\HH$
 containing $v_i$ and deleting $v_i$ from all of them.
We have $|\HH_i|=d_i$.
The hypergraph $\HH_i$ contains at least $d_i - \ex(\KK_{d-1}^{(k)},K_{(d-1)\ast 2})$ copies of $K_{(d-1)\ast 2}$.
Since $\HH$ is $K_{d\ast 2}$-free, each copy of $K_{(d-1)\ast 2}$ belongs to no more than one of the hypergraphs
 $\HH_1,\HH_2,\dots,\HH_{k^{2^{d-1}}}$.
This implies
\[
 \sum_{i=1}^{k^{2^{d-1}}} \left[ d_i -\left( 2-\frac{1}{2^{d-2}} \right) k^{2^{d-1} -2} \right]
\leq \binom{k}{2} \binom{k^2}{2}\dots \binom{k^{2^{d-2}}}{2}
< \frac{k^{2(2^{d-1}-1)}}{2^{d-1}},
\]
and the claim follows by rearranging the inequality.
\end{proof}

\begin{proof}[Proof of the upper bound in {\rm \bf \tref{bdf}}.]

For  $m=k^{2^d -1}$ we define a family $\FF$ of size $m$ such that
  every subfamily $\FF'$ avoiding $B_d$ has size at most  $2k^{2^d-2}$.
Then $f(m, B_d\text{\rm -free}) \le  O(m^{e_d'})$ follows for all $m$ by the monotonicity of $f$.

Let $\FF$ be a product of $d$ chains, the $i$th of which has size $k^{2^{i-1}}$, i.e., for $1 \le i \le d, 1 \le j \le k^{2^{i-1}}$,
 let $S^i_j$ be sets satisfying\\
$\bullet\ |S^i_j|=j$, $S^i_{j_1} \subset S^i_{j_2}$ if $j_1 \le j_2$,\\
$\bullet\ S^i_{k^{2^{i-1}}}\cap S^j_{k^{2^{j-1}}}=\emptyset$ if $i \ne j$, and\\
$\bullet\ \cF:=\{S^1_{j_1}\cup S^2_{j_2} \cup \dots \cup S^d_{j_d} : 1 \le i \le d, 1 \le j_i \le k^{2^{i-1}}\}$.

Each set in $\FF$ corresponds to a hyperedge in $\KK_d^{(k)}$, and each copy of $B_d$ in $\FF$ corresponds to a copy
 of $\KK_{d\ast 2}$ in $\KK_d^{(k)}$.
The $B_d$-free subfamilies of $\FF$ correspond to $\KK_{d\ast 2}$-free subhypergraphs of $\KK_d^{(k)}$.
The bound  in Theorem~\ref{K22free} on the size of a $\KK_{d\ast 2}$-free subfamily completes the proof.
\end{proof}

\section{Union-free subfamilies}
\label{sec:Uf}

\begin{proof}[Proof of {\rm \bf \tref{afree}}.]
The proof of our lower bound is based on Kleitman~\cite{Kl},
the proof by Erd\H os and Shelah \cite{ES1972} does not work in the general $a$-union free setting.

Let $\FF$ be an arbitrary family of size $m$ and  let $\ell$ be the size of a longest chain in it.
Split $\FF$ according the rank of the sets,  $\FF=\cup_{1\leq k\leq \ell}  \FF_k$.
Each $\FF_k$ together with a chain of size $k$ with a top member from $\FF_k$  form an $a$-union free
 subfamily implying $f(\FF,a\mbox{-union free})\geq |\FF_k|+k-1$ for all $k$.
Adding up we have $\ell\times  f\geq m+ \binom{\ell}{2}$ implying
  $f(\FF,a\mbox{-union free})\geq |\FF|/\ell +(\ell-1)/2$.
Since the lower bound by Fox,  Lee, and Sudakov~\cite{FLS} supersedes ours,   we omit the details.

For the proof of the upper bound (\ref{eq:3}), first we consider the family $\FF_{ES}(k)$ of size $k^2$, what
 Erd\H{o}s and Shelah \cite{ES1972} used to obtain the upper bound (\ref{ES_upper}) on $f(k^2, 2\text{-union free})$.
The family $\FF_{ES}$  is a product of two vertex disjoint chains of lengths $k$, that is,
 given the chains $\emptyset\neq A_1 \subset A_2 \subset \dots \subset A_{k}$ and
 $\emptyset\neq B_1 \subset B_2 \subset \dots \subset B_{k}$ with
 $A_{k} \cap B_{k}=\emptyset$ we define $\FF_{ES}(k):=\{A_i \cup B_j: 1 \le i,j \le k\}$.
We have  $|\FF_{ES}|= k^2$.

\begin{lemma}
\label{lem:ES} If $\cG$ is an $a$-union free subfamily of $\cF_{ES}(k)$, then
    $$|\cG|\le 2(\lceil\sqrt{a+1}\rceil-1)k. $$
\end{lemma}

\begin{proof}  Associate a point set $P$  of the $2$-dimensional grid to the family $\GG$ as
  $P:=\{ (i,j): $ when $A_i \cup B_j\in \GG\}$.
The rectangle $R(i,j)$ is defined as $R(i,j):=\{ (x,y): 1\leq x\leq i$ and $1\leq y\leq j\}$.
The set $A_i\cup B_j$ is a union of $a$ distinct members of $\GG$ if and only if the rectangle $R=R(i,j)$
 contains at least $a$ distinct points apart from $(i,j)$ and
 at least one of these lies on the top boundary of $R$, i.e., on the segment $[(1,j),(i,j)]$ and at least
 one on the rightmost column  $[(i,1),(i,j)]$.

Construct $P'\subseteq P$ by deleting the bottom
$\lceil \sqrt{a+1}\rceil -1$ elements of $P$ in each column of the grid.
Suppose that $P'$ has a row with at least $\lceil \sqrt{a+1}\rceil$ elements, and let $(i,j)$ be the rightmost point.
Then $P$ has at least $\lceil\sqrt{a+1}\rceil^2 \geq a+1$ points in the rectangle  $R(i,j)$, also
 points on the top and the right most sides, a contradiction.
Therefore, $P$ has at most $2(\lceil\sqrt{a+1}\rceil-1)k$ elements.
\end{proof}

Now we are ready to define a family $\FF$ of size $q k^2$, such that
\begin{equation}\label{eq:4}
  f(\FF, a\text{-union free}) <  a-2+2k(\lceil\sqrt{a+1}\rceil-1)+(2k-1)(q-1).
  \end{equation}
The family $\FF$ consists of $q$ \textit{levels}, each of them isomorphic to $\FF_{ES}(k)$.
For all $1 \le \ell \le q$, let
 $\emptyset \neq A^\ell_1 \subset A^\ell_2 \subset \dots \subset A^\ell_{k}$ and
 $\emptyset\neq B^\ell_1 \subset B^\ell_2 \subset \dots \subset B^\ell_{k}$
 be chains of length $k$ such that the $2q$ top sets $A^\ell_{k}$ and $B^{\ell'}_{k}$ are pairwise disjoint.
Let us define
\[
\cF_\ell=\left\{\bigcup_{s=1}^{\ell-1}(A^s_{k} \cup B^s_{k})\cup A^\ell_i \cup B^\ell_j: 1 \le i,j \le k\right\}
\text{ and }
\cF:=\bigcup_{\ell=1}^q\cF_\ell.
\]

Observe that  $|\cF|=m=qk^2$ and indeed each $\cF_\ell$ is isomorphic to $\cF_{ES}$.
Note that if $\ell<\ell'$  and $F \in \cF_\ell, F' \in F_{\ell'}$ then $F \subset F'$.
Let $\cG$ be an $a$-union free subfamily of $\cF$ and let us write $\cG_\ell=\cG \cap \cF_\ell$.
Let $t$ be the smallest integer with $\sum_{\ell=1}^t|\cG_\ell| \ge a-2$.
If there exists no such $t$, then $|\cG| < a-2$, and we are done.
We have:\\
$\bullet$ $\sum_{\ell=1}^{t-1}|\cG_\ell| <a-2$,  by the definition of $t$, \\
$\bullet$ $|\cG_t| \le 2(\lceil\sqrt{a+1}\rceil-1)k$  by \lref{ES} since $\cF_t$ is isomorphic to $\cF_{ES}$, \\
$\bullet$ the family $\GG_\ell$ is 2-union free  for each $\ell$ with $t<\ell\le k$.

To see the latest statement, suppose, on the contrary, that  $G' \cup G''=G$ for some $G,G', G'' \in \cG_\ell$.
Pick any $a-2$ sets $G_1,G_2,\dots,G_{a-2}$ from $\cup_{s=1}^t\cG_s$,  and we have
 $G=G'\cup G''\cup G_1 \cup \dots \cup G_{a-2}$, contradicting $\cG$ being $a$-union free.
Therefore  $|\cG_\ell|\le 2k-1$ by a slight strenghtening of the result of Erd\H os and Shelah
 (see \cite{FLS}).
Putting these observations together, using $|\GG|=\sum|\GG_\ell|$ and $t\geq 1$, we obtain
 (\ref{eq:4}).
Finally,  substituting $q=\lceil \sqrt{a+1}\rceil $ and $k=\lceil \sqrt{m/q}\rceil $ into (\ref{eq:4})
 we have $f(m, a\text{-union free})\leq a+ (4k-1)(2q-1)$.
A little calculation yields (\ref{eq:3}).
\end{proof}

\section{Problems, concluding remarks}
\label{sec:conc}

\begin{conjecture}
If $m=2^n$, then the family consisting of $m$ sets that contains the highest number of subfamilies
 forming a Boolean algebra of dimension $d$ is $2^{[n]}$.
	\end{conjecture}

In Theorem~\ref{K22free} we have considered  $d$-partite hypergraphs
 with very uneven part sizes.  There is a number of results of this type,
 see, e.g.,  Gy\H ori~\cite{Gy}.
Also the sizes grow exponentially, one can easily generalize it for other sequences.

Concerning $a$-union free families we had the modest conjecture
\begin{equation}\label{prob_limit}
  \lim_{a \rightarrow \infty} \left( \liminf_{m \rightarrow \infty}\frac{f(m,a\text{\rm -union~free})}{\sqrt m} \right) \rightarrow \infty
 \end{equation}
Knowing the results of  Fox, Lee, and Sudakov~\cite{FLS} it is natural to ask
\begin{prob}\label{new_limit}
Given $a$, what is the limit
$$  \lim_{m \rightarrow \infty}\frac{f(m,a\text{\rm -union~free})}{a^{1/4}\sqrt m} ?
 $$
\end{prob}
If it exists, it is between $1/3$ and $4$.

One can improve the coefficient $4$ of the factor $a^{1/4}$ in Theorem~\ref{thm:afree}
if in Section~\ref{sec:Uf} we use different sizes. Namely we construct
 $\FF$ by using $\FF_\ell=\FF_{ES}(k_\ell)$ where $k_\ell=k\left(\frac{b-1}{b-2}\right)^{2(\ell-1)}$ with $b=\lceil \sqrt{a+1} \rceil$.
If $q/b$ tends to infinity, we obtain
  $$ f(m, a\text{-union free})\leq \sqrt{8} a^{1/4}\sqrt{m}+O(a).$$

A family $\cF$ is  \textit{$(a,b)$-union free} if there are no distinct sets $F_1, F_2\dots, F_{a+b}$
 satisfying $F_1\cup F_2 \cup \dots \cup F_a=F_{a+1}\cup \dots \cup F_{a+b}$.
This is another frequently investigated property.
However $f(m,(a,b)\text{\rm -free})=a+b-1$ if $a,b \ge 2$, as it is
shown by the family consisting of all $(m-1)$-subsets of an $m$-set.

Many more problems remained open.

\section*{Acknowledgement}
We are greatly indebted to the organizers of  the  $1^{st}$ Eml\'ekt\'abla Workshop July 26--29, 2010,
 Gy\"on\-gy\"os\-tarj\'an, Hungary, where most of the research presented in this paper was done as a
group work.

\end{document}